\documentclass[10pt]{article}
\usepackage{amssymb,amsmath,color}

\allowdisplaybreaks
\begin{document}

\def\pd#1#2{\frac{\partial#1}{\partial#2}}
\def\dfrac{\displaystyle\frac}
\let\oldsection\section
\renewcommand\section{\setcounter{equation}{0}\oldsection}
\renewcommand\thesection{\arabic{section}}
\renewcommand\theequation{\thesection.\arabic{equation}}

\def\Xint#1{\mathchoice
  {\XXint\displaystyle\textstyle{#1}}%
  {\XXint\textstyle\scriptstyle{#1}}%
  {\XXint\scriptstyle\scriptscriptstyle{#1}}%
  {\XXint\scriptscriptstyle\scriptscriptstyle{#1}}%
  \!\int}
\def\XXint#1#2#3{{\setbox0=\hbox{$#1{#2#3}{\int}$}
  \vcenter{\hbox{$#2#3$}}\kern-.5\wd0}}
\def\ddashint{\Xint=}
\def\dashint{\Xint-}

\newcommand{\cred}{\color{red}}
\newcommand{\be}{\begin{equation}\label}
\newcommand{\ee}{\end{equation}}
\newcommand{\bea}{\begin{eqnarray}\label}
\newcommand{\eea}{\end{eqnarray}}
\newcommand{\nn}{\nonumber}
\newcommand{\intO}{\int_\Omega}
\newcommand{\Om}{\Omega}
\newcommand{\cd}{\cdot}
\newcommand{\pa}{\partial}
\newcommand{\ep}{\varepsilon}
\newcommand{\uep}{u_{\eta}}
\newcommand{\vep}{v_{\eta}}
\newcommand{\Sep}{S_{\eta}}
\newcommand{\Lp}{L^p(\Om)}
\newcommand{\Lq}{L^q(\Om)}
\newcommand{\abs}{\\[2mm]}
\newcommand{\eps}{\varepsilon}
\newtheorem{thm}{\indent Theorem}
\newtheorem{lem}{\indent Lemma}[section]
\newtheorem{proposition}{\indent Proposition}[section]
\newtheorem{dnt}{\indent Definition}[section]
\newtheorem{remark}{\indent Remark}[section]
\newtheorem{cor}{\indent Corollary}[section]
\allowdisplaybreaks

\title{Global bounded solutions of the higher-dimensional Keller-Segel system under smallness conditions in optimal spaces}
\author{
Xinru Cao\footnote{caoxinru@gmail.com}\\
{\small School of Mathematical Science,
Dalian University of Technology,}\\
{\small 116023 Dalian, China} \\
{\small and}\\
{\small Institut f\"ur Mathematik, Universit\"at Paderborn,}\\
{\small 33098 Paderborn, Germany}}
{\small }
{\small }
\maketitle
\date{}
{\bf Abstract} In this paper, the fully parabolic Keller-Segel system
\be{a0}
\left\{
\begin{array}{llc}
u_t=\Delta u-\nabla\cdot(u\nabla v),
&(x,t)\in \Omega\times (0,T),\\
v_t=\Delta v-v+u, &(x,t)\in\Omega\times (0,T),\\
\end{array}
\right.
	\qquad \qquad (\star)
\ee
is considered under Neumann boundary conditions in a bounded
domain $\Om\subset\mathbb{R}^n$ with smooth boundary,
where $n\ge 2$. We derive a smallness condition on the initial data in optimal Lebesgue spaces
which ensure global boundedness and large time convergence.
More precisely, we shall show that one can find $\varepsilon_0>0$ such that for all suitably regular initial data
$(u_0,v_0)$ satisfying $\|u_0\|_{L^{\frac{n}{2}}(\Om)}<\varepsilon_0$ and $\|\nabla v_0\|_{L^{n}(\Om)}<\varepsilon_0$,
the above problem possesses a global classical solution which is bounded and converges to the constant steady state
$(m,m)$ with $m:=\frac{1}{|\Om|}\intO u_0$. \\
Our approach allows us to furthermore study a general chemotaxis system with rotational sensitivity in dimension 2,
which is lacking the natural energy structure associated with ($\star$).
For such systems, we prove a global existence and boundedness result under corresponding smallness conditions
on the initially present total mass of cells and the chemical gradient.\abs

\section{Introduction}

In this paper, we consider the initial-boundary value problem for two coupled parabolic equations,
\be{a}
\left\{
\begin{array}{llc}
u_t=\Delta u-\nabla\cdot(u S(u,v,x)\cdot\nabla v), &(x,t)\in \Omega\times (0,T),\\[6pt]
\displaystyle
v_t=\Delta v-v+u, &(x,t)\in\Omega\times (0,T),\\[6pt]
\displaystyle
\nabla v\cdot\nu=(\nabla u-uS(x,u,v)\cdot\nabla v)\cdot\nu=0,&(x,t)\in \partial\Omega\times (0,T),\\[6pt]
u(x,0)=u_{0}(x),\quad  v(x,0)=v_{0}(x),
& x\in\Omega,
\end{array}
\right.
\ee
where $\Omega\subset\mathbb{R}^n$ ($n\ge 2$) is
 a bounded domain
with smooth boundary and $\nu$ is the normal outer vector on $\pa\Om$. $S$ is either a scalar number or ${S}(u,v,x)=({s}_{i,j}(u,v,x))_{n\times n}$
is supposed to be a matrix with
${s}_{i,j}\in
C^2([0,\infty)\times [0,\infty)\times\bar{\Omega})$ $(i,j=1,2)$ and satisfies
\bea{S}
|S(u,v,x)|\le C_S\text{ with }C_S>0.
\eea

Systems of this type describe the evolution of cell populations and their movement affected by the gradients of a chemical signal produced by the cells themselves, a mechanism commonly called chemotaxis.
A classical chemotaxis system was proposed by Keller and Segel \cite{Keller-Segel}. In the simplest form of their model,
the first equation in (\ref{a}) reads
\bea{01}
u_t=\Delta u-\nabla\cdot(\chi u\nabla v),
\eea
where $u$, $v$ denote the density of the cell population and chemical substance concentration, respectively.
The number $\chi \in\mathbb{R}$ measures the sensitivity of the chemotactic response to the chemical gradients, and
the second term on the right of (\ref{01}) reflects the hypothesis that cells move towards {\em higher} densities
of the signal. The second equation in (\ref{a}) models the assumptions that the chemical substance is produced by cells and degrades.
This kind of model has been widely studied during the last 40 years \cite{herrero_velazquez,HW}. We also refer to
the survey \cite{H,HP} for a broad overview.

Among the large quantity of the related researches, deciding whether solutions exist globally or blow up in finite time
seems to be one of the most challenging mathematical topics \cite{Wbl}.
In the two-dimensional setting, a critical mass phenomenon has been identified and studied in many works. In the case
$\intO u_0<4\pi$, the solution is global and bounded \cite{Nagai-Senba-Yoshida}, whereas if $\intO u_0>4\pi$,
the occurrence of blow-up for some initial data is only detected when the second equation is replaced by an
elliptic equation of the form
$-\Delta v+v-u=0$ or $-\Delta v-(u-\bar{u}_0)=0$, which reflects a certain limit procedure \cite{JL}.
In higher dimensions, there are many results for such simplified parabolic-elliptic versions.
For instance, in \cite{C}, it is proved that the corresponding Cauchy problem possesses a
global weak  solution whenever $\|u_0\|_{L^\frac{n}{2}(\mathbb{R}^n)}<C$ with some $C>0$.
Considering the fully parabolic case, it is also shown that for small initial data,
\bea{02}
\|u_0\|_{L^q(\mathbb{R}^n)}<\eps, \text{ and } \|\nabla v_0\|_{L^p(\mathbb{R}^n)}<\eps
\eea
with some suitably small constant $\eps>0$, $q> \frac{n}{2}$ and $p\ge n$,
the solution exists globally and is bounded \cite{CP}.
However, for the fully parabolic version in bounded domains, the same conclusion is up to now known to hold only for
$q>\frac{n}{2}$ and
$p>n$ \cite{W4}. It is our goal to extend this result in the corresponding critical case, that is, for
$q=\frac{n}{2}$ and $p=n$. 	

In contrast to (\ref{01}), a recent study suggests a more general model which allows a wider direction of the cells' movement, such as they move not to the higher density of chemical any more but with a rotation. Then in this system a
sensitivity {\em tensor} $S(u,v,x)$, instead of a {\em scalar} constant $\chi$, is introduced to describe chemotactic motion
\cite{Xue-Othmer}. The introduction of this {\em tensor} valued sensitivity is caused by a kind of complicated interactions between the cell motion speed and directional effects stemming from the action of gravity, for example.

In our study of this new model, we concentrate on the two-dimensional case, and we anticipate that small mass of $u$ guarantees global existence, which indeed parallels the case of a scalar sensitivity.
However, the classical way of proof in the scalar case strongly depends on the use of an energy inequality
\cite{Nagai-Senba-Yoshida}, which is apparently lacking in the general system.
To the best of our knowledge, the only results on global existence and boundedness in a related case
can be found in \cite{Li-Suen-Xue-Winkler}, but with the second equation being replaced by $v_t=\Delta v-f(u)v$,
whereby $v$ enjoys an a priori upper bound according to the obvious estimate $v \le \|v_0\|_{L^\infty(\Om)}$.
Under mild assumptions on $S$ and $f$, the authors in that work proved that the solution exists globally and is
bounded if either $\|v_0\|_{L^\infty(\Om)}$ or $\|u_0\|_{L^1(\Om)}$ is small enough.

Since the second equation of (\ref{a}) has a production term, the method in \cite{Li-Suen-Xue-Winkler}
does not apply to the present situation. However, we may benefit from our approach developed above for
the case of scalar sensitivities in order to prove global boundedness under the assumptions that (\ref{S}) holds,
and that both $\|u_0\|_{L^1(\Om)}$ and $\|\nabla v_0\|_{L^2(\Om)}$ are small enough.
We underline that this assumption on the initial data is still stronger than that in the case of scalar sensitivity,
but our results include exponential convergence, which has not been found before
without assuming $\|\nabla v_0\|_{L^2(\Om)}$ small enough.

Our main result says that
\begin{itemize}
\item if $n\ge 2$, one can find an upper bound for $u_0$ in $L^\frac{n}{2}(\Om)$ and $\nabla v_0$ in $L^n(\Om)$,
such that the solution $(u,v)$ of \eqref{a} exists globally and is bounded.
\end{itemize}

\section{Preliminaries}
In this section, we recall some classical $L^p-L^q$ estimates for the
Neumann heat semigroup on bounded domains. Almost all of the results and their proofs can be found in \cite[Lemma 1.3]{W4}. However, some of the estimates we use below go slightly beyond,
and since we could not find a precise reference, we will give a short proof here.

\begin{lem}
Suppose $(e^{t\Delta})_{t>0}$ is the Neumann heat semigroup
in $\Om$, and let $\lambda_1>0$ denote the first nonzero
 eigenvalue of $-\Delta$ in $\Om$ under Neumann
boundary conditions. Then there exist $k_1,...,k_4>0$ which only
depend on $\Om$ and which have the following properties:

\begin{itemize}
\item[\rm(i)]
If $1\le q\le p\le\infty$, then
\bea{1}
\|e^{t\Delta}w\|_{\Lp}\le k_1(1+t^{-\frac{n}{2}
(\frac{1}{q}-\frac{1}{p})})e^{-\lambda_1 t}\|w\|_{\Lq}
\text{ for all }t>0
\eea
holds for all $w\in\Lq$ with $\intO w=0$.
\item[\rm(ii)]
If $1\le q\le p\le\infty$, then
\bea{2}
\|\nabla e^{t\Delta}w\|_{\Lp}\le k_2(1+t^{-\frac{1}{2}-\frac{n}{2}
(\frac{1}{q}-\frac{1}{p})})e^{-\lambda_1 t}\|w\|_{\Lq}
\text{ for all }t>0
\eea
holds for each $w\in\Lq$.
\item[\rm(iii)]
If $2\le q\le p<\infty$, then
\bea{3}
\|\nabla e^{t\Delta}w\|_{\Lp}\le k_3e^{-\lambda_1 t}(1+t^{-\frac{n}{2}(\frac{1}{q}-\frac{1}{p})})
\|\nabla w\|_{\Lq}
\text{ for all }t>0
\eea
is true for all $w\in W^{1,p}(\Om)$.
\item[\rm(iv)]
Let $1<q\le p\le\infty$, then
\bea{4}
\|e^{t\Delta}\nabla\cdot w\|_{\Lp}\le k_4(1+t^{-\frac{1}{2}
-\frac{n}{2}(\frac{1}{q}-\frac{1}{p})})e^{-\lambda_1 t}\|w\|_{\Lq}
\text{ for all }t>0
\eea
is valid for any $w\in (W^{1,p}(\Om))^n$.
\end{itemize}
\end{lem}
{\bf Proof}~~\rm(i) and \rm(ii)  are precisely proved in \cite[Lemma 1.3]{W4}. Focusing on \rm(iii), we note that it is obviously true for all $t<2$ \cite{W4}. If $t\ge2$, let $\bar{w}=\frac{1}{|\Om|}\intO w$, we use \eqref{1}, \eqref{2} and the Poinc\'are inequality to obtain
\begin{align}\nn
\|\nabla e^{t\Delta}w\|_{\Lp}&=\|\nabla e^{\Delta}e^{(t-1)\Delta}(w-\bar{w})\|_{\Lp}\\\nn
&\le 2k_2\|e^{(t-1)\Delta}(w-\bar{w})\|_{\Lp}\\\nn
&\le 2k_2\cdot k_1(1+(t-1)^{-\frac{n}{2}
(\frac{1}{q}-\frac{1}{p})})e^{-\lambda_1(t-1)}
\|w-\bar{w}\|_{\Lq}\\\nn
&\le k_3(1+t^{-\frac{n}{2}
(\frac{1}{q}-\frac{1}{p})})e^{-\lambda_1t}
\|\nabla w\|_{\Lq}.
\end{align}
Thus \rm(iii) is valid for all $t>0$. Note that $k_3$ is independent of $p$ or $q$ since the constant in Poinc\'are inequality could be independent of $q$.

In \cite{W4}, (\ref{4}) is proved for $1<q\le p<\infty$, so that it is sufficient to prove the case $1<q<p=\infty$.
Suppose that $w\in (C_0^\infty(\Om))^n$. Then $\intO e^{\frac{t}{2}\Delta}\nabla\cdot w=\intO \nabla\cdot w=0$, whence
from \eqref{1} we see that
\begin{align*}\nn
\|e^{t\Delta}\nabla\cdot w\|_{L^\infty(\Om)}
&=\|e^{\frac{t}{2}\Delta}(e^{\frac{t}{2}\Delta}\nabla\cdot w)\|
_{L^\infty(\Om)}
\le k_1(1+(\frac{t}{2})^{-\frac{n}{2}\frac{1}{q}})
e^{-\lambda_1\frac{t}{2}}\|
e^{\frac{t}{2}\Delta}\nabla\cdot w\|_{L^q(\Om)}\\\nn
&\le k_1(1+(\frac{t}{2})^{-\frac{n}{2}\frac{1}{q}})
e^{-\lambda_1\frac{t}{2}}\cdot
k_2(1+(\frac{t}{2})^{-\frac{1}{2}})e^{-\lambda_1\frac{t}{2}}
\|w\|_{L^q(\Om)}\\\nn
&\le k_3(1+t^{-\frac{1}{2}-\frac{n}{2q}})
e^{-\lambda_1t}\|w\|_{L^q(\Om)}
\end{align*}
for all $t\in(0,T)$ with $k_3=3k_1k_2$. Thus \eqref{4} is obtained by means of a unique extension to all of $(W^{1,p}(\Om))^n$.
\quad$\Box$\medskip

Before going into the main part, we also recall some local existence and extensibility results for \eqref{a}. We refer to \cite[Propersition 2.3]{CS} for the case with rotation. See also \cite[Lemma 1.1]{W4} for the case where $S$ is a scalar constant. Beyond these, we generalize as follows.
\begin{lem}\label{lem}
Suppose that either $S(x,u,v)\in \mathbb{R}$ or  $S(x,u,v)$ is a matrix satisfying \eqref{S} and in addition that $S(x,u,v)=0$ on $\pa\Om$. Assume
$u_0\in C(\bar{\Om})$ and $v_0\in W^{1,\sigma}(\Om)$ are nonegative with $\sigma>n$. Then there exists $T_{\max}>0$ such that \eqref{a} possesses a classical solution $(u,v)\in \big(C^0([0,T_{max})\times\bar{\Omega})
\cap C^{2,1}((0,T_{max})\times\bar{\Omega})\big)^2$.
Moreover, $T_{\max}<\infty$ if and only if
\begin{align}
  \mathop{\limsup}\limits_{{t\nearrow T_{max}}}
  \|u(\cdot,t)\|_{L^{\infty}(\Omega)}=\infty.
\end{align}
\end{lem}

Now we give a sufficient condition for boundedness and global existence.

\begin{lem}\label{bdlem}
Let $\theta>\frac{n}{2}$. If the solution of (\ref{a}) satisfies
\bea{bdd}
\|u\|_{L^\theta(\Om)}< \infty\text{ for all } t\in(0,T_{\max}),
\eea
then $T_{\max}=\infty$, and $\mathop{\sup}\limits_{t>0}(\|u\|_{L^\infty(\Om)}
+\|v\|_{W^{1,\infty}(\Om)})<\infty$. \end{lem}

The following lemma presents an estimate for certain integrals, which we will frequently use in the next section.
The proof can be found in \cite[Lemma 1.2]{W4}.

\begin{lem}\label{inte}
There exists $C>0$ such that for $\alpha<1$, $\beta<1$ and $\gamma$, $\delta$ be positive and satisfy $\gamma\neq\delta$. Then
\bea{ie}
\int_0^t(1+(t-s)^{-\alpha})e^{-\gamma(t-s)}(1+s^{-\beta})
e^{-\delta s}ds\le C(\alpha,\beta,\delta,\gamma)
(1+t^{\min\{0,1-\alpha-\beta\}})e^{-\min\{\gamma,\delta\}\cdot t}
\eea
for all $t>0$, where $C(\alpha,\beta,\delta,\gamma):=2C\cdot(\frac{1}{|\delta-\gamma|}+
\frac{1}{1-\alpha}+\frac{1}{1-\beta})$.
\end{lem}

\section{Smallness conditions in optimal spaces}
Now having at hand the tools collected in the last section, we can prove global
existence in the classical Keller-Segel system
\be{b}
\left\{
\begin{array}{llc}
\displaystyle
u_t=\Delta u-\nabla\cdot(u\nabla v),
&(x,t)\in \Omega\times (0,T),\\
\displaystyle
v_t=\Delta v-v+u, &(x,t)\in\Omega\times (0,T),\\
\displaystyle
\nabla u\cdot \nu=\nabla v\cdot \nu=0, &(x,t)\in\pa\Om\times (0,T),\\
\displaystyle
u(x,0)=u_0(x), v(x,0)=v_0(x),&x\in\Om,
\end{array}
\right.
\ee
where $\Om\subset\mathbb{R}^n$ with smooth boundary and $n\ge2$.

Our proof relies on a fixed-point type argument developed in \cite{W4}. Unlike in
the original proof there, our assumption is that the initial data be suitably small in {\em optimal} spaces
in the sense that we require that
\bea{}\nn
\|u_0\|_{L^{\frac{n}{2}}(\Om)} \text{ and } \|\nabla v_0\|_{L^n(\Om)} \text{ are small.}
\eea
This seems too weak to give an
$L^\infty(\Om)$ bound in a one-step procedure.
However, we can first use
the "weak" assumption to obtain the smallness condition in supercritical spaces, which meet the assumptions
in \cite[Theorem 2.1]{W4}. Finally, we can derive convergence in $L^\infty(\Om)$, and obtain the convergence rate $e^{-\lambda't}$ with any $0<\lambda'<\lambda_1$, where $\lambda_1>0$ denote the first nonzero
 eigenvalue of $-\Delta$ in $\Om$ under Neumann
boundary conditions.

In Theorem \ref{pro1}, we show how we improve the smallness condition into a supercritical space.
\begin{thm}\label{pro1}
Let $n\ge2$, $0<\lambda'<\lambda_1$. Then there exists $\eps_0>0$ depending on $\lambda'$ and $\Om$ with the following property:
If $u_0\in C^0(\bar{\Om})$ and $v_0\in W^{1,\sigma}(\Om)$
with $\sigma>n$ are nonnegative and satisfy
\bea{11}
\|u_0\|_{L^\frac{n}{2}(\Om)}\le\varepsilon \text{     and     }
\|\nabla v_0\|_{L^n(\Om)}\le\varepsilon
\eea
for some $\varepsilon<\varepsilon_0$, then (\ref{b}) possesses a global classical solution $(u,v)$ which
is bounded and satisfies
\bea{1.0}
\|u(\cdot,t)-e^{t\Delta}u_0\|
_{L^\theta(\Om)}\le\varepsilon
(1+t^{-1+\frac{n}{2\theta}})e^{-\lambda't}
\text{ for all } t>0,
\eea
for all $\theta\in[q_0,\theta_0]$, where $\frac{n}{2}<q_0<n$ and
 $n<\theta_0<\frac{nq_0}{n-q_0}$.
\end{thm}

{\bf Proof} ~~First we fix $0<\lambda'<\lambda_1$. Since
$\frac{n}{2}<q_0<n$ and
 $n<\theta_0<\frac{nq_0}{n-q_0}$, it is possible to fix $q_1,q_2>0$ satisfying $q_0< q_1<\theta_0$,
 $n<q_2<\frac{nq_0}{n-q_0}$ and $\frac{1}{q_0}=\frac{1}{q_1}+\frac{1}{q_2}$.
By applying Lemma \ref{inte}, we can find constants $c_1,c_2,c_3>0$ such that for all $\theta\in[q_0,\theta_0]$,
\begin{align}\label{1.1}
&\int_0^t(1+(t-s)^{-\frac{1}{2}-\frac{n}{2}(\frac{1}{\theta}-\frac{1}{p})})
e^{-(\lambda_1+1)(t-s)}(1+s^{-1+\frac{n}{2\theta}})
e^{-\lambda's}ds \le c_1(1+t^{\min\{0,-\frac{1}{2}+\frac{n}{2p}\}})
e^{-\lambda't}\\\label{1.2}
&\int_0^t(1+(t-s)^{-\frac{1}{2}-\frac{n}{2}
(\frac{1}{q_0}-\frac{1}{\theta})})e^{-\lambda_1(t-s)}
(1+s^{-\frac{3}{2}+\frac{n}{2q_0}})
e^{-\lambda's}ds\le c_2(1+t^{\min\{0,-1+\frac{n}{2\theta}\}})e^{-\lambda't}\\
\label{1.3}
&\int_0^t(1+(t-s)^{-\frac{1}{2}-\frac{n}{2}(\frac{1}{q_0}-\frac{1}{\theta})})
e^{-\lambda_1(t-s)}(1+s^{-\frac{1}{2}+\frac{n}{2q_0}})
e^{-\lambda's}ds\le c_3(1+t^{\min\{0,\frac{n}{2\theta}\}})e^{-\lambda't}
\end{align}
hold for all $t>0$. Since $p<\frac{nq_0}{n-q_0}$, we know that $\frac{1}{2}-
\frac{n}{2}(\frac{1}{\theta_0}-\frac{1}{p})
>\frac{n}{2}(\frac{1}{q_0}-\frac{1}{\theta_0})>0$, and $\lambda_1+1-\lambda'>1$, then we see that
$c_1>0$ is dependent on $\Om$ only.
Similarly, it is not difficult to derive that
$\frac{1}{2}-\frac{n}{2}(\frac{1}{q_0}-
\frac{1}{\theta})\ge\frac{1}{2}-\frac{n}{2}
(\frac{1}{q_0}-\frac{1}{\theta_0})>0$ and $1-1+\frac{n}{2q_1}-\frac{1}{2}
+\frac{n}{2q_2}=-\frac{1}{2}+\frac{n}{2q_0}>0$ due to $q_0<n$ and $\theta_0<\frac{nq_0}{n-q_0}$.
Therefore, $c_2>0$ depend only on $\lambda'$, $\Om$. More precisely, $c_2\to\infty$ as $\lambda'\to\lambda_1$.
Noting that
$1-\frac{1}{2}+\frac{n}{2q_0}>\frac{1}{2}>0$, $c_3$ depends on $\lambda'$ and $\Om$.

 With $\eps_0>0$ to be determined below, we assume \eqref{11} holds for $\eps\in(0,\eps_0)$, and define
\bea{12}
T:=\sup\big\{\tilde{T}>0\big| \|u(\cdot,t)-e^{t\Delta}u_0\|
_{L^\theta(\Om)}\le\varepsilon
(1+t^{-1+\frac{n}{2\theta}})e^{-\lambda't},
\text{ for all }t\in[0,\tilde{T}), \text{ for all }\theta\in[q_0,\theta_0]\big\}.
\eea
We observe that $T$ is well defined and positive due to Lemma \ref{lem}. It is sufficient to prove that $T=\infty$.
Since $n\ge 2$, it is easy to see that
\bea{u0}
\bar{u}_0=
\frac{1}{|\Om|}\intO u_0\le |\Om|^{-\frac{2}{n}}\|u_0\|_{L^\frac{n}{2}(\Om)}\le |\Om|^{-\frac{2}{n}}\eps.
\eea
Then by the definition of $T$ and (\ref{1}), for each $\theta\in[q_0,\theta_0]$,
\begin{align}\nn
\|u(\cdot,t)-\bar{u}_0\|_{L^\theta(\Om)}&\le\|u-e^{t\Delta}u_0\|
_{L^\theta(\Om)}+\|e^{t\Delta}u_0-\bar{u}_0\|_{L^\theta(\Om)}\\\nn
&\le\varepsilon(1+t^{-1+\frac{n}{2\theta}})e^{-\lambda't}
+k_1(1+t^{-1+\frac{n}{2\theta}})e^{-\lambda_1t}
\|u_0-\bar{u}_0\|_{L^{\frac{n}{2}}(\Om)}\\\nn
&\le(\varepsilon+k_1\|u_0\|_{L^{\frac{n}{2}}(\Om)}
+k_1\bar{u}_0|\Om|^\frac{2}{n})(1+t^{-1+\frac{n}{2\theta}})e^{-\lambda't}\\\label{utheta}
&\le c_4\varepsilon(1+t^{-1+\frac{n}{2\theta}})e^{-\lambda't}
\end{align}
holds for all $t\in[0,T)$, where $c_4=1+2k_1$ only depends on
$\Om$.

Next, we claim that there exists $c_5>0$ such that for all $p\in[q_0,\frac{nq_0}{n-q_0})$,
\bea{nv}
\|\nabla(v(\cdot,t)-e^{t(\Delta-1)}v_0)\|_{L^{p}(\Om)}\le
c_5\eps(1+t^{-\frac{1}{2}+\frac{n}{2p}})e^{-\lambda't} \text{ for all }
t\in[0,T),\eea
where $c_5>0$ depends only on $\Om$.
Indeed,
by using the variation-of-constants formula for $v$,
\bea{v}\nn
v(\cdot,t)=e^{t(\Delta-1)}v_0+\int_0^t e^{(t-s)(\Delta-1)}u(\cdot,s)ds ~~~~\text{ for all } t\in(0,T)
\eea
and with the fact that $\nabla e^{t\Delta}\bar{u}_0=\nabla\bar{u}_0=0$, we obtain
\begin{align}\nn
\|\nabla(v(\cdot,t)-e^{t(\Delta-1)}v_0)\|_{L^{p}(\Om)}
&\le\int_0^t\|\nabla e^{(t-s)(\Delta-1)}(u(\cdot,s)-\bar{u}_0)\|_{L^{p}(\Om)}ds
\end{align}
for all $t\in[0,T)$.
We observe that (\ref{2}) (\ref{utheta}) and (\ref{1.1}) imply
\begin{align}\nn
&~~\|\nabla(v(\cdot,t)-e^{t(\Delta-1)}v_0)\|_{L^{p}(\Om)}\\\nn
&\le\int_0^t\|\nabla e^{(t-s)(\Delta-1)}(u(\cdot,s)-\bar{u}_0)\|
_{L^{p}(\Om)}ds\\\nn
&\le \int_0^t k_2(1+(t-s)^{-\frac{1}{2}-\frac{n}{2}(
\frac{1}{\theta}-\frac{1}{p})})
e^{-(\lambda_1+1)(t-s)}
\|u-\bar{u}_0\|_{L^{\theta}(\Om)}ds\\\nn
&\le \int_0^t k_2(1+(t-s)^{-\frac{1}{2}-\frac{n}{2}(
\frac{1}{\theta}-\frac{1}{p})})
e^{-(\lambda_1+1)(t-s)}
c_4\eps(1+s^{-1+\frac{n}{2\theta}})
e^{-\lambda's}ds\\\nn
&\le k_2c_4c_1\eps(1+t^{-\frac{1}{2}+\frac{n}{2p}})
e^{-\lambda't}
\end{align}
for all $t\in[0,T)$ and all $p\in[q_0,\frac{nq_0}{n-q_0})$. Now we gain \eqref{nv} by letting $c_5=k_2c_1c_4$, where $c_5>0$ depends on $\Om$ only.

On the other hand, by applying (\ref{3}) and making use of (\ref{11}), we obtain
\begin{align}\label{nv1}
\|\nabla e^{t(\Delta-1)}v_0\|_{L^{p}(\Om)}&\le k_3e^{-(\lambda_1+1)t}(1+t^{-\frac{1}{2}
+\frac{n}{2p}})
\|\nabla v_0\|_{L^n(\Om)}
\le k_3\varepsilon (1+t^{-\frac{1}{2}+\frac{n}{2p}})
e^{-(\lambda_1+1)t}
\end{align}
for all $t\in[0,T)$ and $p\in[n,\frac{nq_0}{n-q_0})$. We observe that a slightly modified version of \eqref{nv1} is also true for $p<n$, because H\"older's inequality, \eqref{2} and \eqref{11} yield
\begin{align}\nn
\|\nabla e^{t(\Delta-1)}v_0\|_{\Lp}&\le  |\Om|^{\frac{1}{p}-\frac{1}{n}}\|\nabla e^{t(\Delta-1)}v_0\|_{L^n(\Om)}\le k_2\max\{1,|\Om|,|\Om|^{-1}\} e^{-(\lambda_1+1)t}\|\nabla v_0\|_{L^n(\Om)}\\\label{nv2}
&\le c_6\varepsilon
e^{-(\lambda_1+1)t}\le c_6\varepsilon(1+t^{-\frac{1}{2}+\frac{n}{2p}})
e^{-\lambda't}
\end{align}
for all $t\in(0,T)$ with $c_6=k_2\max\{1,~|\Om|\}$ independent of $p$.
Thus collecting \eqref{nv}, \eqref{nv1} and \eqref{nv2}, and observing that $\lambda_1+1>\lambda'$,  we conclude that there exists $c_7>0$ such that
\begin{align}\nn
\|\nabla v(\cdot,t)\|_{L^{p}(\Om)}&\le c_5\varepsilon(1+t^{-\frac{1}{2}+
\frac{n}{2p}})e^{-\lambda't}
+(c_6+k_3)\varepsilon (1+t^{-\frac{1}{2}+\frac{n}{2p}})
e^{-(\lambda_1+1)t}\\\label{vnabla}
&\le c_7\varepsilon(1+t^{-\frac{1}{2}+\frac{n}{2p}})
e^{-\lambda't}
\end{align}
for all $t\in[0,T)$ and all $p\in[q_0,\frac{nq_0}{n-q_0})$, where $c_7=c_5+c_6+k_3$ depends on $\Om$ only.

Now we use the variation-of-constants formula associated with $u$ to infer that
\begin{align}\nn
\|u(\cdot,t)-e^{t\Delta}u_0\|_{L^\theta(\Om)}&\le \int_0^t
\|e^{(t-s)\Delta}\nabla\cdot(u\nabla v)\|_{L^\theta(\Om)}ds\\\nn
&\le \int_0^tk_4(1+(t-s)^
{-\frac{1}{2}-\frac{n}{2}
(\frac{1}{q_0}-\frac{1}{\theta})})
e^{-\lambda_1(t-s)}
\|u(\cdot,s)\nabla v(\cdot,s)\|_{L^{q_0}(\Om)}ds\\\nn
&\le \int_0^tk_4(1+(t-s)
^{-\frac{1}{2}-\frac{n}{2}
(\frac{1}{q_0}-\frac{1}{\theta})})
e^{-\lambda_1(t-s)}
\|(u-\bar{u}_0)\nabla v(\cdot,s)\|_{L^{q_0}(\Om)}ds\\\nn
&~~~~~~+
\int_0^tk_4(1+(t-s)
^{-\frac{1}{2}-\frac{n}{2}
(\frac{1}{q_0}-\frac{1}{\theta})})
e^{-\lambda_1(t-s)}\bar{u}_0\|\nabla v(\cdot,s)\|_{L^{q_0}(\Om)}ds\\\label{u}
&=:I_1+I_2
\end{align}
for all $t\in[0,T)$, where $k_4$ from \eqref{4} depends on $\Om$ only. To estimate $I_1$, we recall the choices of $q_1$ and $q_2$, which ensure that $-1+\frac{n}{2q_1}<0$,
$-\frac{1}{2}
+\frac{n}{2q_2}<0$ and $-1+\frac{n}{2q_1}
-\frac{1}{2}
+\frac{n}{2q_2}=-\frac{3}{2}+\frac{n}{2q_0}>-1$ as well as
$-\frac{1}{2}-\frac{n}{2}(\frac{1}{q_0}
-\frac{1}{\theta})>
-\frac{1}{2}-\frac{n}{2}(\frac{1}{q_0}
-\frac{1}{\theta_0})>-1$.
Moreover, since $q_0<q_1<\theta_0$ and $n<q_2<\frac{nq_0}{n-q_0}$,
in light of the H\"{o}lder inequality, \eqref{utheta} \eqref{vnabla} and \eqref{1.2}, we conclude for all $\theta\in[q_0,\theta_0]$ that
\begin{align}\nn
I_1&\le\int_0^tk_4(1+(t-s)^
{-\frac{1}{2}-\frac{n}{2}
(\frac{1}{q_0}-\frac{1}{\theta})})
e^{-\lambda_1(t-s)}
\|u(\cdot,s)-\bar{u}_0\|_{L^{q_1}(\Om)}\|\nabla v(\cdot,s)\|_{L^{q_2}(\Om)}ds\\\nn
&\le \int_0^tk_4(1+(t-s)^
{-\frac{1}{2}-\frac{n}{2}
(\frac{1}{q_0}-\frac{1}{\theta})})
e^{-\lambda_1(t-s)}
c_4\varepsilon(1+s^{-1+\frac{n}{2q_1}})
e^{-\lambda's}
c_7\varepsilon(1+s^{-\frac{1}{2}
+\frac{n}{2q_2}})e^{-\lambda's}
ds\\\nn
&\le \int_0^t4c_4c_7k_4\varepsilon^2
(1+(t-s)^
{-\frac{1}{2}-\frac{n}{2}
(\frac{1}{q_0}-\frac{1}{\theta})})
e^{-\lambda_1(t-s)}
(1+s^{-1+\frac{n}{2q_1}
-\frac{1}{2}
+\frac{n}{2q_2}})e^{-\lambda's}
ds\\\nn
&\le 4c_4c_7k_4c_2\varepsilon^2(1+t^
{\min\{0,1-\frac{1}{2}-\frac{n}{2}
(\frac{1}{q_0}-\frac{1}{\theta})
-1+\frac{n}{2q_1}
-\frac{1}{2}+\frac{n}{2q_2}\}})
e^{-\lambda't}\\\label{I1}
&= c_8\varepsilon^2(1+t^{-1+\frac{n}{2\theta}})
e^{-\lambda't}
\end{align}
for all $t\in[0,T)$ with $c_8=4c_4c_7k_4c_2$. We see that $c_8>0$ depend only on $\lambda'$, $\Om$. More precisely, $c_8\to\infty$ as $\lambda'\to\lambda_1$.

Similarly, \eqref{vnabla} together with \eqref{1.3} imply
\begin{align}\nn
I_2&\le \int_0^t k_4(1+(t-s)^{-\frac{1}{2}-\frac{n}{2}
(\frac{1}{q_0}-\frac{1}{\theta})})
e^{-\lambda_1(t-s)}
\bar{u}_0\|\nabla v(\cdot,s)\|_{L^{q_0}}ds\\\nn
&\le \int_0^tk_4(1+(t-s)^{-\frac{1}{2}
-\frac{n}{2}
(\frac{1}{q_0}-\frac{1}{\theta})})
e^{-\lambda_1(t-s)}|\Om|^{-\frac{2}{n}}
\varepsilon\cdot
c_7\varepsilon(1+s^{-\frac{1}{2}
+\frac{n}{2q_0}})e^{-\lambda's}
ds\\\nn
&\le c_7k_4c_3|\Om|^{-\frac{2}{n}}\varepsilon^2(1+
t^{\min\{0,1-\frac{1}{2}-
\frac{n}{2}(\frac{1}{q_0}-\frac{1}{\theta})
-\frac{1}{2}+\frac{n}{2q_0}\}})
e^{-\lambda't}\\\nn
&\le c_7k_4c_3|\Om|^{-\frac{2}{n}}\varepsilon^2(1+t^{\min\{0,
\frac{n}{2\theta}\}})e^{-\lambda't}\\\label{I2}
&\le 2c_7k_4c_3|\Om|^{-\frac{2}{n}}\varepsilon^2
e^{-\lambda't}
\le c_9\varepsilon^2(1+t^{-1+\frac{n}{2\theta}})
e^{-\lambda't}
\end{align}
for all $t\in[0,T)$ with $c_9=2c_7k_4c_3|\Om|^{-\frac{2}{n}}$.
Thus $c_9$ depends on $\Om$ and $\lambda'$, and $c_9\to\infty$ as $\lambda'\to\lambda_1$.
As a consequence of \eqref{u}, \eqref{I1} and \eqref{I2}, we arrive at
\begin{align}
\|u(\cdot,t)-e^{t\Delta}u_0\|_{L^\theta(\Om)}\le c_{10}\varepsilon^2
(1+t^{-1+\frac{n}{2\theta}})
e^{-\lambda't}\text{ for all }t\in[0,T), \theta\in[q_0,\theta_0]
\end{align}
where $c_{10}=c_8+c_9$ depends only on $\lambda'$ and $\Om$.
Choosing $\varepsilon_0<\frac{1}{2c_{10}}$, we conclude that $T=\infty$.

Since $T\le T_{\rm max}$, we conclude that $T_{\rm max}=\infty$. In order to prove $(u,v)$ is bounded, we first obtain from the definition of $T$ and \eqref{u0} that
\bea{uu}
\|u(\cdot,t)\|_{L^\theta(\Om)}\le \|u(\cdot,t)-\bar{u}_0\|_{L^\theta}
+\bar{u}_0|\Om|^{\frac{1}{\theta}}
\le c_4\eps(1+t^{-1+\frac{n}{2\theta}})
e^{-\lambda't}+|\Om|^{\frac{1}{\theta}-\frac{2}{n}}\eps
\eea
is true for $t\in[1,\infty)$ for any choice of $\theta\in[q_0,\theta_0)$.
The result of local existence yields
\bea{tsmall}
\|u(\cdot,t)\|_{L^\theta(\Om)}<\infty
\eea
for all $t<1$. We see that (\ref{tsmall}) together with \eqref{uu} yields that there exists $c_{11}>0$ such that
\bea{ub}
\|u(\cdot,t)\|_{L^\theta(\Om)}\le c_{11}
 \eea
 for all $\theta\in[q_0,\theta_0)$ and all $t>0$. Since $\theta\ge q_0>\frac{n}{2}$, \eqref{ub}
implies boundedness of $(u,v)$ by Lemma \ref{bdlem}. Thus $(u,v)$ is global and bounded.
\quad$\Box$\medskip

\begin{remark}
A careful re-inspection of the above argument shows that for the constant $c_{10}=c_{10}(\lambda')$ satisfies
$c_{10}(\lambda')\to\infty$ as $\lambda'\to\lambda_1$,
where the constant $\eps_0=\eps_0(\lambda')$ by Theorem \ref{pro1} has the property that $\eps_0(\lambda')\to 0$
as
$\lambda'\to\lambda_1$.
\end{remark}

With the decay property in higher Lebesgue spaces obtained above, we can obtain a smallness condition
which ensures stabilization of $(u,v)$ in $L^\infty(\Om)$.

\begin{thm}\label{th1}
Let $n\ge2$, $0<\lambda'<\lambda_1$, then there exists $\hat{\eps}_0>0$ with the following property:

If $u_0\in C^0(\bar{\Om})$ and $v_0\in W^{1,\sigma}(\Om)$
with $\sigma>n$ are nonnegative and
\bea{20}
\|u_0\|_{L^\frac{n}{2}(\Om)}\le\varepsilon \text{     and     }
\|\nabla v_0\|_{L^n(\Om)}\le\varepsilon
\eea
for some $\varepsilon\le\hat{\varepsilon}_0$, then (\ref{b}) possesses a global classical solution $(u,v)$ which
is bounded and satisfies
\bea{21}
\|u(\cdot,t)-\bar{u}_0\|_{L^\infty(\Om)}\le \bar{C}e^{-\lambda't}\text{ and  }
\|v(\cdot,t)-\bar{u}_0\|_{L^\infty(\Om)}\le \bar{C}e^{-\min\{\lambda',1\}t}\text{ for all } t>0,
\eea
where $\bar{C}>0$ depends on $\lambda'$ and $\Om$.

\end{thm}
{\bf Proof}~~~
First we fix $q_0\in(\frac{n}{2},n)$, $\theta_0\in(n,\frac{nq_0}{n-q_0})$ and $\theta\in[q_0,\theta_0]$ as well as $p\in[q_0,\frac{nq_0}{n-q_0}]$. In particular,
we know that $\theta>\frac{n}{2}$ and $p>n$.
By applying Theorem 2.1 in \cite{W4}, we see that there exists
$\eps_1>0$, such that for any nonnegative initial data $(\tilde{u}_0,\tilde{v}_0)\in C(\bar{\Om})\times W^{1,\sigma}(\Om)$ with $\sigma>n$ satisfying
\bea{til0}
\|\tilde{u}_0\|_{L^\theta(\Om)}\le\eps_1,
\quad\|\nabla\tilde{v}_0\|_{L^p(\Om)}\le\eps_1,
\eea
then the solution $(\tilde{u},\tilde{v})$ of (\ref{b}) with the initial data $(\tilde{u}_0,\tilde{v}_0)$ exists globally and satisfies
\bea{til}\nn
&&\|\tilde{u}(\cdot,t)-e^{t\Delta}{\tilde{u}}_0\|_{L^\infty(\Om)}\le C\eps^2e^{-\lambda' t},\\\label{tiluv}
&&\|\tilde{u}(\cdot,t)-\overline{\tilde{u}}_0\|_{L^\infty}\le C e^{-\lambda't},
\text{ and }\|\tilde{v}(\cdot,t)-\overline{\tilde{u}}_0\|_{L^\infty(\Om)}\le C e^{-\min\{\lambda',1\}t}
\eea
for all $t>0$ and some $C>0$. According to  (\ref{utheta}) and (\ref{vnabla}), there exists $\eps_0>0$ such that for any initial data satisfying (\ref{20}),
\bea{unv}\nn
&&\|u(\cdot,t)\|_{L^\theta(\Om)}\le\|u(\cdot,t)-\bar{u}_0\|_{L^\theta(\Om)}
+\bar{u}_0|\Om|^{\frac{1}{\theta}}\le \eps_0(1+t^{-1+\frac{n}{2\theta}})e^{-\lambda't}
+|\Om|^{\frac{1}{\theta}-\frac{2}{n}}\eps_0,\\\nn
&&\|\nabla v(\cdot,t)\|_{L^p(\Om)}\le c_7\eps_0(1+t^{-1+\frac{n}{2p}})e^{-\lambda't}
\eea
hold for all $t>0$. From this, we observe that there exists $t_0>0$ such that for all $t\ge t_0$,
\bea{unv1}\nn
&&\|u(\cdot,t)\|_{L^\theta(\Om)}\le 2|\Om|^{\frac{1}{\theta}-\frac{2}{n}}\eps_0,\\\nn
&&\|\nabla v(\cdot,t)\|_{L^p(\Om)}\le 2|\Om|^{\frac{1}{\theta}-\frac{2}{n}}\eps_0.
\eea
Let $\tilde{u}_0=u(t_0)$, $\tilde{v}_0=v(t_0)$, we easily see that $\overline{\tilde{u}}_0=\bar{u}_0$, $\tilde{u}(t)=u(t+t_0)$ and $\tilde{v}(t)=v(t+t_0)$. Taking $\hat{\eps}_0=\min\{\eps_0,\frac{1}{2}|\Om|^{\frac{2}{n}-\frac{1}{\theta}}\eps_1\}$ and substituting $(u,v)$ into \eqref{tiluv} will complete the proof.
\quad$\Box$\medskip

\section{System with rotational sensitivity}
In this section, we consider the modified Keller-Segel system with general tensor-valued sensitivity as given by
\be{ba}
\left\{
\begin{array}{llc}
u_t=\Delta u-\nabla\cdot(uS(u,v,x)\cdot\nabla v),
&(x,t)\in \Omega\times (0,T),\\[6pt]
\displaystyle
v_t=\Delta v-v+u,
&(x,t)\in\Omega\times (0,T),\\[6pt]
\displaystyle
(\nabla u-uS(u,v,x)\cdot\nabla v)\cdot\nu=0,\quad \nabla v\cdot\nu=0
&(x,t)\in \partial\Omega\times (0,T),\\[6pt]
u(x,0)=u_{0}(x),\quad  v(x,0)=v_{0}(x),
& x\in\Omega,
\end{array}
\right.
\ee
where $\Om\subset\mathbb{R}^2$ with smooth boundary. The sensitivity $S$
is now supposed to be a tensor-valued function satisfying (\ref{S}).
The non-flux and coupled boundary condition complicate the solvability.
Following \cite{Li-Suen-Xue-Winkler}, we first regularize the system as below
\be{be}
\left\{
\begin{array}{llc}
(\uep)_t=\Delta \uep-
\nabla\cdot(\uep S_\eta(\uep,\vep,x)\cdot\nabla \vep),
&(x,t)\in \Omega\times (0,T),\\[6pt]
(\vep)_t=\Delta \vep-\vep+\uep,
&(x,t)\in\Omega\times (0,T),\\[6pt]
\displaystyle
(\nabla \uep-
\uep S_\eta(\uep,\vep,x)\cdot\nabla\vep)\cdot\nu=0,\quad
\nabla \vep\cdot\nu=0
&(x,t)\in \partial\Omega\times (0,T),\\[6pt]
\uep(x,0)=u_{0}(x),\quad  \vep(x,0)=v_{0}(x),
& x\in\Omega,
\end{array}
\right.
\ee
with $\Sep(x,\uep,\vep)=\rho_\eta(x)S(\uep,\vep,x)$,
which vanishes on the boundary $\pa\Om$ if $\rho_\eta$ is a suitable cut-off function on $\Om$. Where $\rho_\eta\in[0,1]$ and satisfies
\bea{rho}
\rho_\eta\to 1 \text{ a.e. as } \eta\to 0.
\eea
The first boundary condition of (\ref{ba})
is then reduced to $\nabla\uep\cdot\nu=0$, so that
local classical solvability of (\ref{ba}) can be obtained by the standard approach (Lemma\ref{lem}). Upon combining the idea in the previous section with a limiting procedure $\eta\to 0$, we will derive the following.

\begin{thm}\label{th2}
Let $S$ satisfy (\ref{S}) and $0<\lambda'<\lambda_1$. Then there exists $\epsilon>0$ with the property that if the nonnegative initial data $u_0\in C(\bar{\Om})$ and $v_0\in W^{1,\sigma}(\Om)$ satisfy
\bea{30}
\|u_0\|_{L^1(\Om)}\le\varepsilon \text{     and     }
\|\nabla v_0\|_{L^2(\Om)}\le\varepsilon
\eea
then (\ref{ba}) possesses a global classical solution $(u,v)$ which is bounded and satisfies
\bea{uvcon}
\|u-\bar{u}_0\|_{L^\infty(\Om)}\le \bar{C}e^{-\lambda' t},
\|v-\bar{u}_0\|_{L^\infty(\Om)}\le \bar{C}e^{-\min\{\lambda',1 \}t}
\eea
with some $\bar{C}>0$.
\end{thm}

Before we proceed to prove Theorem \ref{th2}, we start with studying the regularized problem (\ref{be}).
Since $|S_\eta(\uep,\vep,x)|\le C_S$, in light of Lemma \ref{lem}, we can apply a slightly modified version of Theorem \ref{th1} to obtain global existence and boundedness for (\ref{be}) under appropriate smallness conditions on the initial data.

\begin{proposition}
Suppose $0<\lambda'<\lambda_1$. Then there exists $\eps_0>0$ with the following property:

If $u_0\in C^0(\bar{\Om})$ and $v_0\in W^{1,\sigma}(\Om)$
with $\sigma>2$ are nonnegative and satisfy
\bea{31}
\|u_0\|_{L^1(\Om)}\le\varepsilon \text{     and     }
\|\nabla v_0\|_{L^2(\Om)}\le\varepsilon
\eea
for some $\varepsilon\le\varepsilon_0$, where $\eps_0$ depends on $\Om$ and $\lambda'$.
Then the classical solution $(\uep,\vep)$ of \eqref{be}
exists globally and stays bounded.
Moreover, there exists $M>0$ depending on $\Om$ and $\lambda'$ such that
\bea{34}& \|\uep-\bar{u}_0\|_{L^\infty(\Om)}\le M\eps e^{-\lambda' t},
\|\vep-\bar{u}_0\|_{L^\infty(\Om)}\le M\eps e^{-\min\{\lambda',1 \}t},\\
\label{bd}
&\|\uep\|_{L^\infty(\Om)}\le M,
\|\vep\|_{L^\infty(\Om)}\le M, \\
\text{ and  }
\label{uv} &\|\uep-\vep\|_{L^\infty(\Om)}\le 2M\eps e^{-\min\{\lambda',1\}t},  \\
\text{as well as }
&\label{35}\|\nabla \vep\|_{L^2(\Om)}\le M\varepsilon e^{-\lambda' t}\text{ for all }t>0.
\eea

\end{proposition}

Note that the above estimates are independent of $\eta$. Having obtained global existence and long time convergence for \eqref{be}, we proceed to construct a
solution of (\ref{ba}) upon letting $\eta\rightarrow 0$. This limit procedure needs some compactness properties of
$(\uep,\vep)$, which are proven in the following lemmata.

\begin{lem}\label{lem1}
There exists $C_1>0$ such that
\begin{align}\label{32}
\int_0^\infty\intO |\nabla\vep|^2\le C_1,\\
\label{33}
\int_0^\infty\intO |\nabla\uep|^2\le C_1.
\end{align}
\end{lem}
{\bf Proof}~~The boundedness of $\int_0^\infty\intO |\nabla\vep|^2$ is an immediate consequence of (\ref{35}). Next, we multiply the first equation of (\ref{be}) by $\uep$ to see that
\begin{align}\nn
\frac{1}{2}\frac{d}{dt}\intO\uep^2+\intO|\nabla\uep|^2
&=\intO\uep S_\eta(\uep,\vep,x)\nabla \vep\cdot\nabla\uep\\
\nn
&\le \frac{1}{2}\intO|\nabla\uep|^2
+\frac{1}{2}C_S^2M^2\intO|\nabla\vep|^2
\end{align}
for all $t>0$.
Rearranging and integrating over $(0,\infty)$ imply
\begin{align}
\int_0^\infty\intO|\nabla\uep|^2
\le\intO u_0^2+C_S^2M^2\int_0^\infty\intO|\nabla\vep|^2.
\end{align}
Finally, we can choose $C_1>0$ in an obvious way to establish \eqref{32} and \eqref{33}.\quad$\Box$\medskip

The next lemma will ensure that $\frac{\pa}{\pa t}\uep$ and $\frac{\pa}{\pa t}\vep$ are bounded in $L^2([0,\infty);(W^{1,2}(\Om))^*)$.

\begin{lem}\label{lem2}
There exists $C_2>0$ such that
\begin{align}\label{4.2.1}
\|\frac{\pa}{\pa t}\uep\|_{L^2([0,\infty);(W^{1,2}(\Om))^*)}\le C_2;\\\label{4.2.2}
\|\frac{\pa}{\pa t}\vep\|_{L^2([0,\infty);(W^{1,2}(\Om))^*)}\le C_2.
\end{align}
\end{lem}
{\bf Proof}~~Let $\phi\in W^{1,2}(\Om)$, and take $\phi$ as a test function in the first equation to obtain from \eqref{35}
\begin{align}\nn
\intO \frac{\pa}{\pa t}\uep\phi
&=-\intO\nabla\uep\cdot\nabla\phi
+\intO\uep S_{\eta}\nabla\vep\cdot\nabla\phi\\
\nn
&\le
\Big(\intO|\nabla\uep|^2\Big)^\frac{1}{2}
\Big(\intO|\nabla\phi|^2\Big)
^\frac{1}{2}
+MC_S\Big(\intO|\nabla\vep|^2\Big)^\frac{1}{2}
\Big(\intO|\nabla\phi|^2\Big)^\frac{1}{2}\\
\nn
&\le \Big(\Big(\intO|\nabla\uep|^2\Big)^\frac{1}{2}
+MC_S\Big(\intO|\nabla\vep|^2\Big)^\frac{1}{2}\Big)
\|\phi\|_{W^{1,2}(\Om)}.
\end{align}
This implies that
$$\|\frac{\pa}{\pa t}\uep(\cdot,t)\|_{(W^{1,2}(\Om))^*}
\le \Big(\intO|\nabla\uep|^2\Big)^\frac{1}{2}
+MC_S\Big(\intO|\nabla\vep|^2\Big)^\frac{1}{2}$$
for all $t>0$, and hence
\bea{ut}\int_0^\infty\|\frac{\pa}{\pa t}\uep(\cdot,t)\|
_{(W^{1,2}(\Om))^*}^2
\le \int_0^\infty\intO|\nabla\uep|^2+M^2C_S^2\int_0^\infty
\intO|\nabla\vep|^2.\eea
The right-hand side of \eqref{ut} is bounded due to Lemma \ref{lem1}.
Similarly, we multiply $\phi$ on both sides of the second
 equation to obtain
\begin{align}\nn
\intO \frac{\pa}{\pa t}\vep\phi
&=-\intO\nabla\vep\cdot\nabla\phi-\intO\vep\phi
+\intO\uep\phi\\\nn
&\le \Big(\intO|\nabla\vep|^2\Big)^\frac{1}{2}
\Big(\intO|\nabla\phi|^2\Big)^\frac{1}{2}
+\intO(\uep-\vep)\phi\\
\nn
&\le \Big(\Big(\intO|\nabla\vep|^2\Big)^\frac{1}{2}
+\Big(\intO|\uep-\vep|^2\Big)^\frac{1}{2}\Big)\|\phi\|_{W^{1,2}(\Om)}.
\end{align}
Again, this fact together with (\ref{34}) entails that
\bea{36}\nn
\|\frac{\pa}{\pa t}\vep(\cdot,t)\|_{(W^{1,2}(\Om))^*}
\le \displaystyle\Big(\intO|\nabla\vep|^2
\Big)^\frac{1}{2}
+\Big(\intO|\uep-\vep|^2\Big)^\frac{1}{2}.
\eea
We integrate over $(0,\infty)$ to obtain
\bea{vt}
\int_0^\infty\|\frac{\pa}{\pa t}\vep(\cdot,t)\|_{(W^{1,2}(\Om))^*}^2
\le \displaystyle\int_0^\infty\intO|\nabla\vep|^2
+\int_0^\infty\intO|\uep-\vep|^2.
\eea
Collecting \eqref{uv}, \eqref{ut}, \eqref{vt}, \eqref{32} and \eqref{33}, we see that if we choose
\bea{37}\nn C_2:=(2+M^2 C_S^2)C_1+\frac{2M^2|\Om|}{\min\{\lambda_1,1\}},
\eea
then \eqref{4.2.1} and \eqref{4.2.2} hold.\quad$\Box$\medskip

Now we can obtain the desired compactness
properties of $(\uep,\vep)$ to prove Theorem \ref{th2}.

{\bf Proof of Theorem \ref{th2}}~~~First we note that both $u_\eta$ and $v_\eta$ are bounded in $L^2_{loc}([0,\infty);W^{1,2}(\Om))$ according to Lemma \ref{lem1} and \eqref{bd}. This fact together with Lemma \ref{lem2} yields that the families $\{u_\eta\}$ and $\{v_\eta\}$  are strongly compact in $L^2_{loc}([0,\infty);L^2(\Om))$ by invoking a version of the Aubin-Lions lemma \cite{Temam}.
We see that there exists $\{\eta_j\}_{j\in \mathbb{N}}\in(0,1)$ satisfying $\eta_j\to 0$ as $j\to\infty$ and nonnegative functions $u,v\in L^2_{loc}([0,\infty);L^2(\Om))$ such that
\bea{41}
u_\eta\to u, \;v_\eta\to v \text{ in }
L^2_{loc}([0,\infty);L^2(\Om)) \text{ as } \eta=\eta_j\to 0.
\eea
According to \eqref{bd}, \eqref{32} and \eqref{33}, we obtain the following properties of $(u,v)$
\bea{38}\nn
&&\uep\stackrel{\star}\rightharpoonup u, \;\vep\stackrel{\star}\rightharpoonup v\text{ in }L^\infty((0,\infty)\times\Om)\\
\label{a.e.}
&&\uep\rightarrow u, \;\vep\rightarrow v\text{ a.e. in } \Om\times(0,\infty),\\\nn
&\text{ and } &\nabla\uep\rightharpoonup\nabla u,\;\nabla\vep\rightharpoonup\nabla v\text{ in }L^2((0,T)\times\Om).
\eea
We see that (\ref{a.e.}) in combination with (\ref{rho}) implies
\bea{ae}
\uep\Sep(\uep,\vep,x)\rightarrow uS(u,v,x) \text{ a.e. in }\Om\times(0,\infty).
\eea
Choosing $\phi\in C_0^\infty(\bar{\Om}\times[0,\infty))$, we see that $(\uep,\vep)$ also satisfies
\bea{ep}
&&-\int_0^\infty\intO\uep\phi_t-\intO u_0\phi(\cdot,0)=-\int_0^\infty\intO\nabla\uep\cdot\nabla\phi
\\\nn
&&~~~~~~~~~~~~~~~~~~~~~~~~~~~~~~~~~~~~~~~~~~~~~~~~+\int_0^\infty\intO \uep(\Sep(\uep,\vep,x)\nabla\vep)\cdot\nabla\phi\\
&&-\int_0^\infty\intO \vep\phi_t-\intO v_0\phi(\cdot,0)=-\int_0^\infty\intO\nabla\vep\cdot\nabla\phi\\\nn
&&~~~~~~~~~~~~~~~~~~~~~~~~~~~~~~~~~~~~~~~~~~~~~~~~-\int_0^\infty\intO \vep\phi+\int_0^\infty\intO\uep\phi.
\eea
Here, (\ref{a.e.}) and (\ref{ae}) allow us to take $\eta=\eta_j\to 0$ in the above identities. Therefore $(u,v)$ is weak solution of (\ref{a}) in the natural sense.
By standard parabolic theory \cite{Ladyenskaja,Li-Suen-Xue-Winkler}, $(u,v)$ is in fact a classical solution of (\ref{ba}). Moreover,
\eqref{a.e.} enable us to take $\eta=\eta_j\to 0$ in \eqref{34} to derive \eqref{uvcon}.\quad$\Box$\medskip

\section{Acknowledgement}
The author would like to thank Professor Michael Winkler for guidance and useful comments on this
paper, which largely improve the recent work.

\end{document}